\documentstyle{article}

\newcommand{\aA}{\mbox{$\cal A$}}

\newcommand{\Ka}{\mbox{$\cal K$}}

\newcommand{\eS}{\mbox{$\cal S$}}

\newcommand{\str}{\rightarrow}
\newcommand{\rts}{\leftarrow}
\newcommand{\mj}{\mbox{\bf 1}}

\newcommand{\df}{\mbox{\scriptsize{\it df}}}
\newcommand{\HDS}{\vrule width0pt height2.3ex depth1.05ex\displaystyle}

\newcommand{\Be}{\mbox{$\cal B$}}

\def\cirk{\,{\raisebox{.3ex}{\tiny $\circ$}}\,}

\def\f#1#2{{{\HDS #1}\over{\HDS #2}}}

\def\nav#1#2{\parshape=2 1em 31em 4em 28em \smallskip
\noindent{\makebox[3em][l]{#1}}{#2}\par\smallskip}

\newcommand{\Jo}{\mbox{$\cal J$}}

\begin{document}

\title{Negation and Involutive Adjunctions}
\author{{\sc Kosta Do\v sen} and {\sc Zoran Petri\' c}\\[0.5cm]
Mathematical Institute, SANU \\
Knez Mihailova 35, p.f. 367 \\
11001 Belgrade, Serbia \\
email: \{kosta, zpetric\}@mi.sanu.ac.yu}
\date{}
\maketitle

\begin{abstract}
\noindent This note analyzes in terms of categorial proof theory
some standard assumptions about negation in the absence of any
other connective. It is shown that the assumptions for an
involutive negation, like classical negation, make a kind of
adjoint situation, which is named involutive adjunction. The
notion of involutive adjunction amounts in a precise sense to
adjunction where an endofunctor is adjoint to itself.
\end{abstract}

\vspace{.3cm}


\vspace{.5ex}

\noindent {\it Keywords$\,$}: negation, adjunction,
self-adjunction

\vspace{0.5cm}

\begin{flushright}
{\it Dedicated to Dov Gabbay on the occasion of his 60th birthday}
\end{flushright}

\vspace{0.5cm}

\baselineskip=1.2\baselineskip

\section{Introduction}

The goal of this note is to present a phenomenon of adjunction
present in assumptions about an involutive negation connective,
like classical negation. Proof-theoretical assumptions concerning
such a negation make an adjoint situation that we call an
\emph{involutive adjunction}. The notion of involutive adjunction
amounts, in a sense to be made precise, to adjunction where an
endofunctor is adjoint to itself, which in \cite{DP03} is called
\emph{self-adjunction}.

In a series of papers, which starts with \cite{Gab88} (see
\cite{Gab91}, \cite{Gab96} and \cite{Gab99}), Dov
Gabbay\index{Gabbay, D. M.} has been working on characterizations
of negation in terms of assumptions about a consequence relation.
Sometimes, as in this note, Gabbay\index{Gabbay, D. M.}
concentrates on negation in the absence of any other connective.
The context of the present note replaces Gabbay's\index{Gabbay, D.
M.} logical framework of a consequence relation by a consequence
\emph{graph},\index{consequence graph} as this is done in
categorial proof theory. We do not have any more only a relation
between premises and conclusions, but we have arrows between them,
and there may be more than one such arrow. We are interested in
equalities between these arrows. Often these equalities, which are
proof-theoretically motivated, exemplify important notions of
category theory. This note shows that with an involutive negation
we fall on a particular notion of adjunction. This is yet another
corroboration of Lawvere's\index{Lawvere, F. W.} thesis that all
logical constants are tied to adjoint situations (see
\cite{Law69}), and of Mac Lane's\index{Mac Lane, S.} slogan that
adjunction arises everywhere (see \cite{ML71}, Preface).

\section{Self-adjunctions}

To fix notation and terminology, we will rely on the following
definition of the notion of adjunction (cf.\ \cite{ML71}, Section
IV.1, and \cite{D99}, Section 4.1.3).

An \emph{adjunction}\index{adjunction} is a sextuple
${\langle\aA,\Be,F,G,\varphi,\gamma\rangle}$ where

\vspace{1ex}

\nav{}{\aA\ and \Be\ are categories,}

\nav{}{$F$ from \Be\ to \aA\ and $G$ from \aA\ to \Be\ are
functors,}

\nav{}{$\varphi$ is a natural transformation of \aA\ from the
composite functor ${FG}$ to the identity functor of \aA, which
means that the following equation holds in \aA\ for every arrow
${f\!:A_1\str A_2}$ of \aA:}

\vspace{-2ex}

\begin{tabbing}
\quad\quad\quad\quad\quad\quad\quad\=($\varphi$~{\it nat})
\quad\quad\= $f\cirk\varphi_{A_1}=\varphi_{A_2}\cirk FGf$,\kill

\>($\varphi$~{\it nat}) \> $f\cirk\varphi_{A_1}=\varphi_{A_2}\cirk
FGf$,

\end{tabbing}

\nav{}{$\gamma$ is a natural transformation of \Be\ from the
identity functor of \Be\ to the composite functor ${GF}$, which
means that the following equation holds in \Be\ for every arrow
${g\!:B_1\str B_2}$ of \Be:}

\vspace{-2ex}

\begin{tabbing}
\quad\quad\quad\quad\quad\quad\quad\=($\varphi$~{\it nat})
\quad\quad\= $f\cirk\varphi_{A_1}=\varphi_{A_2}\cirk FGf$,\kill

\>($\gamma$~{\it nat}) \> $GFg\cirk\gamma_{B_1}=\gamma_{B_2}\cirk
g$,

\end{tabbing}

\nav{}{the following \emph{triangular} equations\index{triangular
equations} hold in \aA\ and \Be\ respectively:}

\vspace{-1ex}

\begin{tabbing}
\quad\quad\quad\quad\quad\quad\quad\=($\varphi$~{\it nat})
\quad\quad\= $f\cirk\varphi_{A_1}=\varphi_{A_2}\cirk FGf$,\kill

\>($\varphi\gamma F$)\>$\varphi_{FB}\cirk
F\gamma_B\:$\=$=\mj_{FB}$,
\\*[1ex]
\>($\varphi\gamma G$)\>$G\varphi_A\cirk\gamma_{GA}$\>$=\mj_{GA}$.

\end{tabbing}

\vspace{1ex}

A \emph{self-adjunction}\index{self-adjunction} is a quadruple
${\langle\eS,L,\varphi,\gamma\rangle}$ where
${\langle\eS,\eS,L,L,\varphi,\gamma\rangle}$ is an adjunction
(this notion is taken over from \cite{DP03}, Section 10). So, in a
self-adjunction, $L$ is an endofunctor, and the equations
($\varphi$~{\it nat}) and ($\gamma$~{\it nat}) become

\begin{tabbing}
\quad\quad\quad\quad\quad\quad\quad\=($\varphi$~{\it nat})
\quad\quad\= $f\cirk\varphi_{A_1}=\varphi_{A_2}\cirk FGf$,\kill

\>\>$f\cirk\varphi_{A_1}=\varphi_{A_2}\cirk LLf$,
\\*[1ex]
\>\>$LLf\cirk\gamma_{A_1}=\gamma_{A_2}\cirk f$,
\end{tabbing}

\noindent while the triangular equations become

\begin{tabbing}
\quad\quad\quad\quad\quad\quad\quad\=($\varphi$~{\it nat})
\quad\quad\= $f\cirk\varphi_{A_1}=\varphi_{A_2}\cirk FGf$,\kill

\>($\varphi\gamma L$)\>$\varphi_{LA}\cirk
L\gamma_A=L\varphi_A\cirk\gamma_{LA}=\mj_{LA}$.

\end{tabbing}

A
\Ka-\emph{self-adjunction}\index{K-self-adjunction@\Ka-self-adjunction}
is a self-adjunction that satisfies the additional equation

\begin{tabbing}
\quad\quad\quad\quad\quad\quad\quad\=($\varphi$~{\it nat})
\quad\quad\= $f\cirk\varphi_{A_1}=\varphi_{A_2}\cirk FGf$,\kill

\>($\varphi\gamma \Ka$)\>$L(\varphi_A\cirk
\gamma_A)=\varphi_{LA}\cirk\gamma_{LA}$,

\end{tabbing}

\noindent and a
\Jo-\emph{self-adjunction}\index{J-self-adjunction@\Jo-self-adjunction}
is a self-adjunction that satisfies the additional equation

\begin{tabbing}
\quad\quad\quad\quad\quad\quad\quad\=($\varphi$~{\it nat})
\quad\quad\= $f\cirk\varphi_{A_1}=\varphi_{A_2}\cirk FGf$,\kill

\>($\varphi\gamma \Jo$)\>$\varphi_A\cirk \gamma_A=\mj_A$

\end{tabbing}

\noindent (these notions are also from \cite{DP03}, Section 10).
It is easy to see that every \Jo-self-adjunction is a
\Ka-self-adjunction (the converse need not hold).

A \Jo-self-adjunction that satisfies

\begin{tabbing}
\quad\quad\quad\quad\quad\quad\quad\=($\varphi$~{\it nat})
\quad\quad\= $f\cirk\varphi_{A_1}=\varphi_{A_2}\cirk FGf$,\kill

\>($\gamma\varphi$)\>$\gamma_A\cirk \varphi_A=\mj_{LLA}$

\end{tabbing}

\noindent is called a \emph{trivial}\index{trivial
self-adjunction} self-adjunction. Note that for trivial
self-adjunctions it is superfluous to assume the equations
($\gamma$~{\it nat}) and ($\varphi\gamma G$), or alternatively
($\varphi$~{\it nat}) and ($\varphi\gamma F$); these equations can
be derived from the remaining ones.

The \emph{free self-adjunction}\index{free self-adjunction}
${\langle\eS,L,\varphi,\gamma\rangle}$ generated by $\{p\}$, where
we call $p$ a \emph{letter}, is defined as follows. The category
\eS\ has as objects the formulae of the propositional language
generated by $\{p\}$ with a unary connective $L$. We may identify
the formulae $p$, $Lp$, $LLp$,$\ldots$ of this language with the
natural numbers 0, 1, 2,$\ldots$

The arrow terms of \eS\ are defined inductively out of the
primitive arrow terms

\[
\mj_A\!:A\str A,\quad\quad \varphi_A\!:LLA\str A, \quad\quad
\gamma_A\!:A\str LLA,
\]

\noindent for every object $A$ of \eS, with the help of the
operations of composition $\cirk$ and the unary operation that
assigns to the arrow term ${f\!:A\str B}$ the arrow term
${Lf\!:LA\str LB}$. On these arrow terms we impose the equations
of self-adjunctions. In the set of these equations we have of
course all the equations $f=f$, and this set is closed under
symmetry and transitivity of equality, and under the rules

\[
\hspace{-6em}(\mbox{\it cong~}\cirk\!)\quad\quad \f{f=f_1 \quad
\quad \quad g=g_1} {f\cirk g=f_1\cirk g_1}
\]

\[
\hspace{-7em}(\mbox{\it cong~L})\quad\quad\f{f=g}{Lf=Lg}
\]

\noindent We assume for $f$ and $g$ in $(\mbox{\it cong~}\cirk\!)$
that they have composable types, such that ${f\cirk g}$ is
defined; the same assumption is made for $f_1$ and $g_1$.

We define analogously the free \Ka-self-adjunction,\index{free
K-self-adjunction@free \Ka-self-adjunction} the free
\Jo-self-adjunc\-tion\index{free J-self-adjunction@free
\Jo-self-adjunction} and the free trivial
self-adjunction\index{free trivial self-adjunction} generated by
$\{p\}$, just by imposing additional equations.

\section{Involutive adjunctions}

Consider a category \aA\ and a contravariant functor
$\neg$\index{negation contravariant endofunctor@$\neg$
contravariant endofunctor} from \aA\ to \aA, which means that for
${f\!:A\str B}$ in \aA\ we have ${\neg f\!:\neg B\str \neg A}$ in
\aA, and the following equations are satisfied:

\begin{tabbing}
\mbox{\hspace{2em}}\= $\mbox{($\neg$)}$\quad\quad\quad\=\kill

\> $\mbox{($\neg 1$)}$\> $\neg\mj_A=\mj_{\neg A}$
\\*[2ex]
\> $\mbox{($\neg 2$)}$\> $\neg(f\cirk g)=\neg g\cirk\neg f$, \quad
for ${f\!:A\str B}$ and ${g\!:C\str A}$.
\end{tabbing}

\noindent The contravariant functor $\neg$ may be conceived either
as a functor from the category $\aA^{op}$ to \aA, which we denote
by $\neg$ too, or as a functor from \aA\ to $\aA^{op}$, which we
denote by $\neg^{op}$.

Suppose that for every object $A$ of \aA\ we have an arrow
${n^{\str}_A\!:\neg\neg A\str A}$ of \aA. The arrow $n^{\str}_A$
becomes the arrow ${n^{\str\:op}_A\!:A\str\neg\neg A}$ in
$\aA^{op}$.

We say that ${\langle\aA,\neg,n^{\str}\rangle}$ is an
$n^{\str}$\emph{-adjunction} when

\[
\langle\aA,\aA^{op},\neg,\neg^{op},n^{\str},n^{\str\:op}\rangle
\]

\noindent is an adjunction. This means that in \aA\ we have for
every ${f\!:A_1\str A_2}$ the equation

\begin{tabbing}
\quad\quad\quad\quad\=(($n^{\str}$~{\it
triang})\quad\quad\=$n^{\str}_{A_1}=n^{\str}_{A_2}\cirk\neg\neg
f$\kill

\>($n^{\str}$~{\it nat})\>$f\cirk
n^{\str}_{A_1}=n^{\str}_{A_2}\cirk\neg\neg f$,

\end{tabbing}

\noindent alternatively written ${f\cirk
n^{\str}_{A_1}=n^{\str}_{A_2}\cirk\neg\neg^{op} f}$, which also
delivers ($n^{\str\:op}$~{\it nat}) in $\aA^{op}$, and the
equation

\begin{tabbing}
\quad\quad\quad\quad\=($n^{\str}$~{\it
triang})\quad\quad\=$n^{\str}_{A_1}=n^{\str}_{A_2}\cirk\neg\neg
f$\kill

\>($n^{\str}$~{\it triang})\>$\;\:n^{\str}_{\neg A}\cirk\neg
n^{\str}_A=\mj_{\neg A}$,

\end{tabbing}

\noindent which delivers both the equation ($\varphi\gamma F$),
i.e.\ ($n^{\str}n^{\str\:op}\;\neg$), in \aA, and the equation
($\varphi\gamma G$), i.e.\ ($n^{\str}n^{\str\:op}\;\neg^{op}$), in
$\aA^{op}$.

Suppose now that we have as before a category \aA\ and a
contravariant functor $\neg$ from \aA\ to \aA, and that for every
object $A$ of \aA\ we have an arrow ${n^{\rts}_A\!:A\str \neg\neg
A}$ of \aA. The arrow $n^{\rts}_A$ becomes the arrow
${n^{\rts\:op}_A\!:\neg\neg A\str A}$ in $\aA^{op}$.

We say that ${\langle\aA,\neg,n^{\rts}\rangle}$ is an
$n^{\rts}$\emph{-adjunction} when

\[
\langle\aA^{op},\aA,\neg^{op},\neg,n^{\rts\:op},n^{\rts}\rangle
\]

\noindent is an adjunction. This means that in \aA\ we have for
every ${f\!:A_1\str A_2}$ the equation

\begin{tabbing}
\quad\quad\quad\quad\=(($n^{\str}$~{\it
triang})\quad\quad\=$n^{\str}_{A_1}=n^{\str}_{A_2}\cirk\neg\neg
f$\kill

\>($n^{\rts}$~{\it nat})\>$\neg\neg f\cirk
n^{\rts}_{A_1}=n^{\rts}_{A_2}\cirk f$,

\end{tabbing}

\noindent which also delivers ($n^{\rts\:op}$~{\it nat}) in
$\aA^{op}$, and the equation

\begin{tabbing}
\quad\quad\quad\quad\=($n^{\str}$~{\it triang})\quad\quad\=$f\cirk
n^{\str}_{A_1}=n^{\str}_{A_2}\cirk\neg\neg f$\kill

\>($n^{\rts}$~{\it triang})\>$\neg n^{\rts}_A\cirk n^{\rts}_{\neg
A}=\mj_{\neg A}$,

\end{tabbing}

\noindent which delivers both the equation ($\varphi\gamma F$),
i.e.\ ($n^{\rts\:op}n^{\rts}\,\neg^{op}$), in $\aA^{op}$, and the
equation ($\varphi\gamma G$), i.e.\
($n^{\rts\:op}n^{\rts}\,\neg$), in \aA. Note that what we call
$n^{\rts}$-adjunction is called \emph{self-adjunction} in
\cite{P97} (Section 3.1; cf.\ also \cite{MLM92}, Section I.8),
which should not be confused with our notion of self-adjunction in
the preceding section.

We say that ${\langle\aA,\neg,n^{\str},n^{\rts}\rangle}$ is an
\emph{involutive}\index{involutive adjunction} adjunction when
${\langle\aA,\neg,n^{\str}\rangle}$ is an $n^{\str}$-adjunction
and ${\langle\aA,\neg,n^{\rts}\rangle}$ is an
$n^{\rts}$-adjunction.

A \Ka\emph{-involutive}\index{K-involutive
adjunction@\Ka-involutive adjunction} adjunction is an involutive
adjunction that satisfies the additional equation

\begin{tabbing}
\quad\quad\quad\quad\=($n^{\str}$~{\it triang})\quad\quad\=$f\cirk
n^{\str}_{A_1}=n^{\str}_{A_2}\cirk\neg\neg f$\kill

\>($n^{\str}n^{\rts}\Ka$)\>$\neg (n^{\str}_A\cirk
n^{\rts}_A)=n^{\str}_{\neg A}\cirk n^{\rts}_{\neg A}$,

\end{tabbing}

\noindent and a \Jo\emph{-involutive}\index{J-involutive
adjunction@\Jo-involutive adjunction} adjunction is an involutive
adjunction that satisfies the additional equation

\begin{tabbing}
\quad\quad\quad\quad\=($n^{\str}$~{\it triang})\quad\quad\=$f\cirk
n^{\str}_{A_1}=n^{\str}_{A_2}\cirk\neg\neg f$\kill

\>($n^{\str}n^{\rts}\Jo$)\>$n^{\str}_A\cirk n^{\rts}_A=\mj_A$.

\end{tabbing}

\noindent It is easy to see that every \Jo-involutive adjunction
is a \Ka-involutive adjunction (the converse need not hold).

A \Jo-involutive adjunction that satisfies

\begin{tabbing}
\quad\quad\quad\quad\=($n^{\str}$~{\it triang})\quad\quad\=$f\cirk
n^{\str}_{A_1}=n^{\str}_{A_2}\cirk\neg\neg f$\kill

\>($n^{\rts}n^{\str}$)\>$n^{\rts}_A\cirk n^{\str}_A=\mj_{\neg\neg
A}$

\end{tabbing}

\noindent is called a \emph{trivial}\index{trivial involutive
adjunction} involutive adjunction.

Note that for trivial  involutive adjunctions it is superfluous to
assume the equations ($n^{\rts}$~{\it nat}) and ($n^{\rts}$~{\it
triang}), or alternatively ($n^{\str}$~{\it nat}) and
($n^{\str}$~{\it triang}); these equations can be derived from the
remaining ones. In trivial involutive adjunctions we have the
equations

\begin{tabbing}
\quad\quad\quad\quad\quad\quad\=($n^{\str}$~{\it
triang})\quad\quad\=$f\cirk
n^{\str}_{A_1}=n^{\str}_{A_2}\cirk\neg\neg f$\kill

\>\>$n^{\rts}_{\neg A}=\neg n^{\str}_A$,
\\*[1ex]
\>\>$n^{\str}_{\neg A}=\neg n^{\rts}_A$.

\end{tabbing}

The \emph{free involutive adjunction}\index{free involutive
adjunction} ${\langle\aA,\neg,n^{\str},n^{\rts}\rangle}$ generated
by $\{p\}$ is defined as follows. The category \aA\ has as objects
the formulae of the propositional language generated by $\{p\}$
with a unary connective $\neg$. We may identify these formulae
with the natural numbers.

The arrow terms of \aA\ are defined inductively out of the
primitive arrow terms

\[
\mj_A\!:A\str A,\quad\quad n^{\str}_A\!:\neg\neg A\str A,
\quad\quad n^{\rts}_A\!:A\str \neg\neg A,
\]

\noindent for every object $A$ of \aA, with the help of the
operations of composition $\cirk$ and the unary operation that
assigns to the arrow term ${f\!:A\str B}$ the arrow term ${\neg
f\!:\neg B\str \neg A}$. On these arrow terms we impose the
equations of involutive adjunctions. In the set of these equations
we have of course all the equations $f=f$, and this set is closed
under symmetry and transitivity of equality, under the rule ({\it
cong}~$\cirk\!$), and also under the rule

\[
\hspace{-7em}(\mbox{\it cong~$\neg$})\quad\quad\f{f=g}{\neg f=\neg
g}
\]

We define analogously the free \Ka-involutive
adjunction,\index{free K-involutive adjunction@free \Ka-involutive
adjunction} the free \Jo-involutive adjunction\index{free
J-involutive adjunction@free \Jo-involutive adjunction} and the
free trivial involutive adjunction\index{free trivial involutive
adjunction} generated by $\{p\}$, just by imposing additional
equations.

Note that the category of the free involutive adjunction generated
by an arbitrary set having more than one letter would be the
disjoint union of isomorphic copies of the category \aA\ of the
free involutive adjunction generated by $\{p\}$. An analogous
remark applies to the category of the free self-adjunction
generated by an arbitrary set having more than one member: it
would be the disjoint union of isomorphic copies of the category
\eS\ of the free self-adjunction generated by $\{p\}$.

\section{Self-adjunctions and
involutive adjunctions}

We are now going to prove that in the free self-adjunction
${\langle\eS,L,\varphi,\gamma\rangle}$ and the free involutive
adjunction ${\langle\aA,\neg,n^{\str},n^{\rts}\rangle}$, both
generated by $\{p\}$, the categories \eS\ and \aA\ are isomorphic
categories.

First, we define $\neg$, $n^{\str}$ and $n^{\rts}$ in \eS\ in the
following manner. On objects we have that $\neg$ is $L$, while for
the arrow term ${f\!:A\str B}$ of \eS\ we define the arrow term
${\neg f\!:\neg B\str \neg A}$ of \eS\ inductively as follows:

\begin{tabbing}
\quad\quad\quad\quad\quad\quad\quad\quad\quad\quad\quad\quad\=$\neg
\varphi_A\:$\=$=L\gamma_A$\kill

\>$\neg\mj_A$\>$=L\mj_A=\mj_{LA}=\mj_{\neg A}$,
\\*[.5ex]
\>$\neg \varphi_A$\>$=L\gamma_A$,
\\*[.5ex]
\>$\neg \gamma_A$\>$=L\varphi_A$,
\\[1.5ex]
\>$\neg(f\cirk g)=\neg g\cirk\neg f$,
\\[.5ex]
\>$\neg Lf$\>$=L\neg f$.

\end{tabbing}

\noindent That this defines an operation $\neg$ on the arrows of
\eS\ is shown by verifying that if ${f=g}$ in \eS, then ${\neg
f=\neg g}$ in \eS; we verify, namely, that the equations of \eS\
are closed under the rule ({\it cong}~$\neg$) of the preceding
section. This is done by a straightforward induction on the length
of the derivation of ${f=g}$ in \eS. For that we use the fact that
for every arrow term $f$ of \eS\ the arrow term $\neg f$ is equal
in \eS\ to an arrow term of the form $Lf'$.

Finally, we have

\[
n^{\str}_A=_{\df}\;\varphi_A,\quad\quad\quad\quad
n^{\rts}_A=_{\df}\;\gamma_A.
\]

Next, we define $L$, $\varphi$ and $\gamma$ in \aA\ in the
following manner. On objects we have that $L$ is $\neg$, while for
the arrow term ${f\!:A\str B}$ of \aA\ we define the arrow term
${L f\!:LA\str LB}$ of \aA\ inductively as follows:

\begin{tabbing}
\quad\quad\quad\quad\quad\quad\quad\quad\quad\quad\quad\quad\=$\neg
\varphi_A\:$\=$=L\gamma_A$\kill

\>$L\mj_A$\>$=\neg\mj_A=\mj_{\neg A}=\mj_{LA}$,
\\*[.5ex]
\>$Ln^{\str}_A$\>$=\neg n^{\rts}_A$,
\\*[.5ex]
\>$Ln^{\rts}_A$\>$=\neg n^{\str}_A$,
\\[1.5ex]
\>$L(f\cirk g)=Lf\cirk Lg$,
\\[.5ex]
\>$L\neg f$\>$=\neg Lf$.

\end{tabbing}

\noindent That this defines an operation $L$ on the arrows of \aA\
is shown by verifying that if ${f=g}$ in \aA, then ${Lf=Lg}$ in
\aA; we verify, namely, that the equations of \aA\ are closed
under the rule ({\it cong}~$L$) of \S 2 above. This is done by a
straightforward induction on the length of the derivation of
${f=g}$ in \aA. For that we use the fact that for every arrow term
$f$ of \aA\ the arrow term $Lf$ is equal in \aA\ to an arrow term
of the form $\neg f'$.

Finally, we have

\[
\varphi_A=_{\df}\;n^{\str}_A,\quad\quad\quad\quad
\gamma_A=_{\df}\;n^{\rts}_A.
\]

We verify easily by induction on the complexity of the arrow term
$f$ that both in \eS\ and in \aA\ we have the equation

\begin{tabbing}
\quad\quad\quad\quad\quad\quad\=\quad\quad\quad\quad\quad\quad\=$\neg
\varphi_A\:$\kill

\>($LL\neg\neg$)\>$LLf=\neg\neg f$.

\end{tabbing}

Next we verify that the equations of involutive adjunctions hold
for the defined $\neg$, $n^{\str}$ and $n^{\rts}$ in \eS. This is
done in a straightforward manner by induction on the length of
derivation. In the basis of this induction, we use ($LL\neg\neg$),
($\varphi$~{\it nat}) and ($\gamma$~{\it nat}) to verify
($n^{\str}$~{\it nat}) and ($n^{\rts}$~{\it nat}), while the
equations ($n^{\str}$~{\it triang}) and ($n^{\rts}$~{\it triang})
reduce to ($\varphi\gamma L$). In the induction step, we rely on
the closure of \eS\ under ({\it cong}~$\neg$), which we
established above.

We verify also that the equations of self-adjunctions hold for the
defined $L$, $\varphi$ and $\gamma$ in \aA. This is done again in
a straightforward manner by induction on the length of derivation.
In the basis of this induction, we use ($LL\neg\neg$),
($n^{\str}$~{\it nat}) and ($n^{\rts}$~{\it nat}) to verify
($\varphi$~{\it nat}) and ($\gamma$~{\it nat}), while the
equations ($\varphi\gamma L$) reduce to ($n^{\str}$~{\it triang})
and ($n^{\rts}$~{\it triang}). In the induction step, we rely on
the closure of \aA\ under ({\it cong}~$L$), which we established
above.

We have a functor $F_{\cal A}$ from \eS\ to \aA\ that maps the
object of \eS\ corresponding to the natural number $n$ to the
object of \aA\ corresponding to $n$, and that maps every arrow of
\eS\ to the homonymous arrow in the defined \eS\ structure of \aA.
For example,

\[
F_{\cal A}\,\varphi_{LLp}=\varphi_{\neg\neg p}=n^{\str}_{\neg\neg
p}\,.
\]

\noindent We define analogously a functor $F_{\cal S}$ from \aA\
to \eS. That $F_{\cal A}$ and $F_{\cal S}$ are indeed functors
follows from what we established above.

It is trivial that on objects we have that ${F_{\cal S}F_{\cal
A}A}$ is $A$, and that ${F_{\cal A}F_{\cal S}B}$ is $B$. We show
next by induction on the complexity of $f$ that in \eS\ we have

\[
F_{\cal S}F_{\cal A}\,f=f.
\]

\noindent When $f$ is of the form $Lf'$, we make an auxiliary
induction on the complexity of $f'$, in which we use
($LL\neg\neg$). We show analogously that in \aA\ we have

\[
F_{\cal A}F_{\cal S}\,g=g.
\]

\noindent This concludes the proof that \eS\ and \aA\ are
isomorphic categories.

We demonstrate analogously that the categories of, respectively,

\nav{}{the free \Ka-self-adjunction and the free \Ka-involutive
adjunction,}

\nav{}{the free \Jo-self-adjunction and the free \Jo-involutive
adjunction,}

\nav{}{the free  trivial self-adjunction and the free
 trivial involutive adjunction,}

 \noindent all generated by $\{p\}$, are isomorphic categories.

 The interest of considering \Ka\ and \Jo\ versions of
 self-adjunctions and involutive adjunctions comes from
 connections with Temperley-Lieb algebras\index{Temperley-Lieb algebras} and the associated
 geometrical interpretation (see \cite{DP03} and references
 therein). Roughly speaking, \Ka\ is what we find in
 Temperley-Lieb algebras,\index{Temperley-Lieb algebras} where only the number of circles (which
 correspond to ${\varphi_A\cirk\gamma_A}$ or ${n^{\str}_A\cirk
 n^{\rts}_A}$) counts, while in \Jo\ circles are disregarded.

 The free trivial self-adjunction, and hence also the free trivial
 involutive adjunction, are preorders; namely, all arrows with the
 same source and target are equal. This follows from the results of
\cite{DP03} (unabridged version) or \cite{DP03Publ}.


\begin{thebibliography}{99}


\bibitem{D99} {\sc K. Do\v sen}, \textbf {\textit
{Cut Elimination in Categories}}, Kluwer, Dordrecht, 1999

\bibitem{DP03} {\sc K. Do\v sen} and {\sc Z. Petri\' c}, {\it Self-adjunctions and
matrices},  \textbf {\textit {Journal of Pure and Applied
Algebra}}, vol.\ 184 (2003), pp.\ 7-39 (unabridged version
available at: http:// arXiv.\-org/\-math.\-GT/\-0111058)

\bibitem{DP03Publ} --------, {\it The geometry of self-adjunction},
\textbf {\textit {Publications de l'Institut Math\' ematique
(N.S.)}}, vol.\ 73(87) (2003), pp.\ 1-29 (incorporated in the
unabridged version of \cite{DP03})

\bibitem{Gab88} {\sc D.~M. Gabbay},
{\it What is negation in a system?}, \textbf {\textit {Logic
Colloquium '86}} (F.~R. Drake and J.~K. Truss, editors),
North-Holland, Amsterdam, 1988, pp.\ 95-112

\bibitem{Gab91} --------,
{\it Modal provability foundations for negation by failure},
\textbf {\textit {Extensions of Logic Programming}} (P.
Schroeder-Heister, editor), Springer, Berlin, 1991, pp.\ 179-222

\bibitem{Gab99} {\sc D.~M. Gabbay} and {\sc A. Hunter},
{\it Negation and contradiction}, \textbf {\textit {What is
Negation?}} (D.M. Gabbay and H. Wansing, editors), Kluwer,
Dordrecht, 1999, pp.\ 89-100

\bibitem{Gab96} {\sc D.~M. Gabbay} and {\sc H. Wansing},
{\it What is negation in a system? Part II: Negation in structured
consequence relations}, \textbf {\textit {Logic, Action and
Information}} (A. Fuhrmann and H. Rott, editors), de Gruyter,
Berlin, 1996, pp.\ 328-350

\bibitem{Law69}  {\sc F.~W. Lawvere}, {\it Adjointness in foundations},
\textbf {\textit {Dialectica}}, vol. 23 (1969), pp. 281-296

\bibitem{ML71} {\sc S. Mac Lane}, \textbf {\textit {Categories for the Working
Mathematician}}, Spring\-er, Berlin, 1971 (expanded second
edition, 1998)

\bibitem{MLM92} {\sc S. Mac Lane} and {\sc I. Moerdijk}, \textbf {\textit
{Sheaves in Geometry and Logic: A First Introduction to Topos
Theory}}, Springer, Berlin, 1992

\bibitem{P97} {\sc D. Pavlovi\' c}, {\it Chu I: cofree equivalences, dualities
and $*$-autonomous categories}, \textbf {\textit {Mathematical
Structures in Computer Science}}, vol.\ 7 (1997), pp.\ 49-73

\end{thebibliography}
\end{document}